\documentclass[12pt,reqno]{amsart}

\title[Wiener-Hopf operators]{Spectra of the translations and Wiener-Hopf operators on $\lw$}

\address{LMAM, Universit\'e de Lorraine (Metz), UMR 7122,Ile du Saulcy 57045, Metz Cedex 1, France.}

\email{petkova@univ-metz.fr}

\usepackage{amsmath} 

\usepackage{amssymb}

\def\squarebox#1{\hbox to #1{\hfill\vbox to #1{\vfill}}}

\newcommand{\F}{{\mathcal F}}

\newcommand{\LL}{{L_{loc}^1(\rr)}}

\newcommand{\I}{[-\alpha_1, \:\alpha_0]}

\newcommand{\rp}{r_+}
\newcommand{\rpb}{{\bf r_+}}
\newcommand{\Vt}{L_{-t}}
\newcommand{\Vu}{L_{-1}}
\newcommand{\St}{R}
\newcommand{\Rd}{{\bf R}}
\newcommand{\Ld}{{\bf L}}

\newcommand{\M}{{\mathcal M}}

\newcommand{\C}{{\mathbb C}}

\newcommand{\N}{{\mathbb N}}

\newcommand{\w}{{\omega}}

\newcommand{\R}{{\Bbb R}}

\newcommand{\lw}{{L_\omega^2({\Bbb R}^+)}}

\newcommand{\rr}{\mathbb R^+}

\newcommand{\Lw}{{L_\omega^2({\Bbb R})}}

\newcommand{\Cc}{{C_c^\infty(\R)}}

\newcommand{\cc}{C_c^\infty(\R^+)}

\newcommand{\CC}{C_c^\infty(\R)}

\renewcommand{\Re}{\mathop{\rm Re}\nolimits}

\renewcommand{\Im}{\mathop{\rm Im}\nolimits}

\theoremstyle{plain}

\def\ccp{C_c(\R^+)}

\newtheorem{thm}{Theorem}

\newtheorem{lem}{Lemma}

\newtheorem{prop}{Proposition}

\newtheorem{deff}{Definition}

\theoremstyle{definition}

\vspace{1cm}

\thanks{Universite de Lorraine, LMAM, UMR 7122, Bat A, Ile de Saulcy, 57045 Metz Cedex 1, France, petkova@univ-metz.fr}
\numberwithin{equation}{section}
\author[Violeta Petkova]{Violeta Petkova}
\begin{document}
\maketitle

\begin{abstract} We study bounded operators $T$ on the weighted space $L^2_{\omega}(\R^+)$ commuting either with the
 "right shift operators" $(\St _t)_{t \geq 0},$ or "left shift operators" $(\Vt)_{t \geq 0},$ and we establish the existence 
of a symbol $\mu$ of $T$. 
We characterize completely the spectrum $\sigma (\St_t)$ of the operator $\St_t$ proving that 
$$\sigma (\St _t) = \{z \in \C: |z| \leq  e^{t\alpha_0}\},$$
 where
 $\alpha_0$ is the growth bound of $(\St_t)_{t\geq 0}$. We obtain a similar result for the spectrum of $(\Vt),\: t >0.$
 Moreover, for a bounded operator $T$ commuting with $\St _t, \: t \geq 0,$ we establish the inclusion 
$\overline{\mu({\mathcal O})}\subset \sigma(T)$, where $\mathcal{O}= \{ z \in \C: \Im z < \alpha_0\}$.

\end{abstract}



\section{Introduction}
Let $\w$ be a weight on $\R^+$. It means that $\w$ is a positive, not vanishing, continuous function on $\R^+$ such that
$$0<\inf_{x\geq 0}\frac{\w(x+t)}{\w(x)}\leq\sup_{x\geq 0}\frac{\w(x+t)}{\w(x)}<+\infty, \forall t\in \R^+.$$
For example $e^x$ and $e^{-x}$ are weights on $\rr$  and we will see later that 
$$\sup_{x\geq 0}\frac{\w(x+t)}{\w(x)}\leq Ce^{mt},$$
where $C$ and $m$ are constants. \\
Let $\lw$ be the set of measurable functions on $\R^+$ such that
$$\int_0^\infty |f(x)|^2\w(x)^2 dx<+\infty.$$
The space $H = \lw$ equipped with the scalar product 
$$<f,g>=\int_{\R^+} f(x)\overline{g}(x)\omega(x)^2dx,\:f\in \lw,\:g\in \lw$$
and the related norm $\|.\|$ is a Hilbert space.
Let $\Cc$ (resp. $\cc$) be the space of $C^\infty$ functions on $\R$ (resp. $\R^+$) with compact support in $\R$
 (resp. $\R^+$). Notice that $\cc$ is dense in $\lw$.
Let $\rp $ be the restriction from $L_{loc}^1(\R^{-})\oplus \lw$ into $\lw$ and let $\rpb $ be the restriction from $L^2(\R)$ into $L^2(\rr)$. 
Let $\ell_0$ be the extension of a function defined on $\rr$ into a function on $\R$ setting
$\ell_0 f(x)=0,\forall x\in \R^-$. These notations are compatible with those used for example in  \cite{Ek}.\\

For $t \in \R$ and a function defined on $\R$ we denote $(S_t f)(x) = f(x - t).$ 
Then for $t > 0$ we introduce the "right shift" operator on $ \lw$ by $\St_t:=\rp S_t \ell_0$.
 For simplicity $\St _1$ will be denoted by $\St $. 
For $t > 0$ we define the "left shift" operator on $\lw$ by $\Vt f:= \rp  S_{-t}\ell_0 f.$ 
It is clear that $\St_t$ and $\Vt$ map $\lw$ into $\lw$, while the shift $S_t$ acts on functions defined on $\R.$
The notations $R_t$ and $L_{-t}$ are similar to those for the operators on the spaces of sequences (see \cite{R1}).
 Let $I$ be the identity operator on $\lw$.


\begin{deff}
A bounded  operator $T$ on $\lw$ is called a Wiener-Hopf operator if 
$$\Vt T\St_tf=Tf,\:\forall t\in \R^+,\:f\in \lw.$$
\end{deff}

More general Wiener-Hopf operators have been intensively studied in the literature (see \cite{Sp} and the references given there).
 There exist also many results about Toeplitz operators on weighted Hardy 
spaces (see \cite{BS}, \cite{BSS}). Such operators have some similarities with Wiener-Hopf operators. \\

Every Wiener-Hopf operator $T$ has a representation by a convolution. The reader may find a proof in \cite{V2} where the arguments of
\cite{H}, \cite{L} are used. More precisely, there exists a distribution $\mu\in\mathcal{D}^\prime(\R)$ such that 
$$Tf(x)=(\mu*\ell_0f)(x),\:\forall f\in C_c^\infty(\rr), \:x\in \rr.$$
If $\phi\in \Cc$ then the operator 
$$\lw\ni f\longrightarrow \rp (\phi*\ell_0 f)$$
is a Wiener-Hopf operator and we will denote it by $T_\phi$.

A bounded operator $T$ on $\lw$ commuting either with $\St _t,\: \forall t > 0$, or with $\Vt,\: \forall t > 0$ is a Wiener-Hopf operator. On the other hand,
every operator $\alpha \Vt + \beta \St _t$ with $t > 0, \:\alpha, \: \beta \in \C$ is a Wiener-Hopf operator. 
It is clear that the set of Wiener-Hopf operators is not a sub-algebra of the algebra of the bounded operators on $\lw$.

Notice also that 
$$(\Vt \St _t)f=f,\:\forall f\in \lw,\: t > 0,$$
but it is obvious that $(\St _t\Vt)f\neq f,$
for all $f\in \lw$ with a support not included in $]t,+\infty[$. 
The fact that $\St_t $ is not invertible leads to many difficulties in contrast to the case when we deal with the space
 $L^2_{\omega}(\R).$ The later space has been considered in \cite{V9} and \cite{V8} and the author has studied  
the operators commuting with the translations on $\Lw$  characterizing their spectrum. 
The group of translations on $\Lw$ is commutative and the investigation of its spectrum is easier. 
In this work, first we apply some ideas used in \cite{V9} and \cite{V8} to study Wiener-Hopf operators on $L_\w^2(\R^+)$. 
For this purpose it is necessary to treat  two semigroups of not invertible operators instead of a group of invertible operators.
 More precisely, we must deal with the semigroups 
$(\St_t)_{t \geq 0}$ and $(\Vt)_{t \geq 0}$ on $L^2_{\omega}(\R^+).$ 
For our analysis it is more convenient to replace the weight $\w$ by another one. 
To do this, given a weight $\w$, denote by $\overline{\w}$ the function 
$\overline{\w}(t)=\sup_{x\geq 0}\frac{\w(x+t)}{\w(x)}.$
Define $\w_0(x)=\exp \Bigl(\int_{1}^{2}\ln( \w(x+t))dt \Bigr).$ Following \cite{BM}, $\w_0$ is a well defined weight on $\rr$ and, moreover,
we have
 \begin{equation}\label{eq:poi}\overline{\w_0}(t)\leq e^{mt},\:\forall t\geq 0,
  \end{equation}
where $m$ is a constant. 
The weights $\w$ and $\w_0$ are equivalent (see \cite{BM}), that is there exist constants $C_1$ and $C_2$ such that 
$$C_1 \w_0(x)\leq \w(x)\leq C_2\w_0(x),\:\forall x\in \rr.$$
Taking into account that $\|\St_t\|=\overline{\w}(t)$, from (\ref{eq:poi}) we get 
 the estimate 
$\|\St_t\|\leq Ce^{m t},\:\forall t\in \rr.$
A similar estimate holds for the semigroup $(\Vt)_{t \geq 0}.$ 

Denote by $\rho(B)$ (resp. $\sigma(B)$) the spectral radius (resp. the spectrum) of an operator $B$.
 Introduce the ground orders
of the semigroups $(\St_t)_{t\geq 0}$ and $(\Vt)_{t\geq 0}$ by
$$\alpha_0 = \lim_{t \to \infty} \frac{1}{t} \ln\|\St_t\|,\:  \alpha_1 = \lim_{t \to \infty} \frac{1}{t}\ln \| \Vt\|.$$
Then it is well known (see for example \cite{EN}) that we have
$$\rho(\St_t) = e^{\alpha_0 t}, \: \rho(\Vt) = e^{\alpha_1 t}.$$
Let $\mathcal{I}$ by the interval $[-\alpha_1, \alpha_0]$ and define
$$\Omega:=\Big\{z\in \C:\: e^{-\alpha_1}\leq |z|\leq e^{\alpha_0}\Big\}.$$
Notice that
$\alpha_1 + \alpha_0 \geq 0$. Indeed, for every $n \in \N$ we have
$L_{-n}\St_{n}=I$ and 
$$1 \leq \limsup_{n \to \infty} \|(\Vu)^n\|^{1/n} \limsup_{n \to \infty} \|\St^n\|^{1/n} = e^{\alpha_1} e^{\alpha_0}$$
which yields the result. For a function $f$ and $a\in \C$ we denote by $(f)_a$ the function 
$$(f)_a: x\longrightarrow f(x)e^{ax}.$$
Denote by $\F$ the usual Fourier transformation on $L^2(\R)$ and set $\hat{f}=\F f$, for $f\in L^2(\R)$.
 For a function $f \in L^2(\R^+),$ we define $\tilde{f}$ as the Fourier transform of the function $f$ extended by 
 0 on $\R^{-}.$ 

Our first result is the following 

\begin{thm}
Let $a\in {\mathcal I} =\I$ and let $T$ be a Wiener-Hopf operator.
There exists $h_a\in L^\infty(\R)$ such that for every $f\in \lw$ satisfying $(f)_a\in L^2(\rr)$, \\

\begin{equation} \label{eq:r}
(Tf)_a=\rpb  \F^{-1}(h_a\widetilde{(f)_a})
\end{equation}
and
$$\|h_a\|_\infty\leq C\|T\|,$$
where $C$ is a constant independent of $a$. Moreover, if $\alpha_1 + \alpha_0 > 0$, the function $h$ defined on $U=\{z\in \C:\:\Im z\in \I\}$ 
by $h(z)=h_{\Im z}(\Re z)$ is holomorphic on $\overset{\circ}{U}$.
\end{thm}
\begin{deff}
The function $h$ defined in Theorem $1$ is called the symbol of $T$. 
\end{deff}

A weaker result that Theorem 1 has been proved in \cite{V2} where the representation (\ref{eq:r}) 
has been obtained only for functions $f \in C_c^\infty(\R^+)$ which is too restrictive for
 the applications to the spectral problems studied in Section 3 and Section 4.
 On the other hand, in the proof in \cite{V2} there is a gap in the approximation argument. 
Indeed in the proof of Lemma 2 in \cite{V2}, there is one argument from Lemma 6 
in \cite{V} which can be applied only if the function $\R^+\ni t\longrightarrow \St_t$
is continuous with respect to the operator norm topology on the set of bounded operators on $\lw$. 
Guided by the approach in \cite{V8}, in this work we prove a stronger version of the result of \cite{V2} 
applying other techniques based essentially on the spectral theory of semigroups.
 On the other hand, in many interesting cases as $\omega(x) = e^x,\: \omega(x) = e^{-x}$, we have $\alpha_0 + \alpha_1 = 0$ 
and the result of Theorem 1 is not satisfied since the symbol of $T$ is defined only on the line $\Im z = \alpha_0.$ To obtain more complete results we introduce the following class of operators.\\
\begin{deff} Denote by $\M$ (resp. $\mathcal{W}$) the set of bounded operators on $\lw$ commuting with 
$\St_t,\: \forall t > 0$ (resp.  $\Vt,\: \forall t > 0$). 
\end{deff}
For operators in $\M$ or $\mathcal{W}$ we obtain a stronger version of Theorem 1.

\begin{thm}
 Let $T$ be a bounded operator commuting with $(\St_t)_ {t > 0}$ (resp. $(\Vt)_{t>0}$). 
Let $a\in J=]0, \alpha_0]$ (resp. $K=]0, \alpha_1]$). 
There exists $h_a\in L^\infty(\R)$ such that for every $f\in \lw$ satisfying $(f)_a\in L^2(\rr)$,
we have 
$$(Tf)_a=\rpb  \F^{-1}(h_a\widetilde{(f)_a})$$
and
$$\|h_a\|_\infty\leq C\|T\|,$$
where $C$ is a constant independent of $a$. Moreover,
 the function $h$ defined on\\
 ${\mathcal O} = \{ z \in \C: \: \Im z < \alpha_0\}$ $($resp. $\mathcal{V}=\{ z \in \C: \: \Im z  >  -\alpha_1\}$\:$)$\\
by $h(z)=h_{\Im z}(\Re z)$ is holomorphic on ${\mathcal O}$ (resp. $\mathcal{V}
$).
\end{thm}

Our main spectral result is the following

\begin{thm} We have
\begin{equation} \label{eq:1.4}
(i)\:\:\:\:\sigma (\St_t) = \{z\in \C,\: |z| \leq e^{\alpha_0 t} \},\:\forall t>0.
\end{equation} 

\begin{equation} \label{eq:1.5}
(ii)\:\:\:\sigma (\Vt) = \{z\in \C,\: |z| \leq e^{\alpha_1 t}\}, \:\forall t>0.
\end{equation}

 iii) Let $T\in \M$ and let $\mu_T$ be the symbol of $T$. Then we have
\begin{equation} \label{eq:1.6}
\overline{\mu_T({\mathcal O})}\subset \sigma(T).
\end{equation}
iv) Let $T\in \mathcal{W}$ and let $\mu_T$ be the symbol of $T$. Then we have
\begin{equation} \label{eq:1.7}
\overline{\mu_T(\mathcal{V})}\subset \sigma(T).
\end{equation}
\end{thm}

It is important to note that for $T \in {\mathcal M}$ (resp. $\mathcal{W}$) and $\lambda\in \C$, 
if the resolvent $(T - \lambda I)^{-1}$ exists, then $(T - \lambda I)^{-1}$ is also in ${\mathcal M}$ (resp. $\mathcal{W}$). In general, if $T$ is a Wiener-Hopf operator and $\lambda\in \C$, even if the resolvent $(T - \lambda I)^{-1}$ exists,
$(T - \lambda I)^{-1}$ could be not a Wiener-Hopf operator. 
For more information about Wiener-Hopf operators the reader may consult \cite{Ek} and \cite{Go}. \\

 The result in Theorem 3 cannot be obtained from a spectral calculus which is unknown and quite 
difficult to construct for the operators in  ${\mathcal M}$ or $\mathcal{W}$. 
On the other hand, our analysis shows the importance of the existence of symbols and this was our main motivation 
to establish Theorem 1 and Theorem 2.\\
The spectrum of the weighted right and left shifts on $l^2(\R^+)$, denoted respectively by $\Rd$ and $\Ld$, has been studied in \cite{R1}. In particular, it was shown that 
\begin{equation} \label{eq:1.3}
\sigma(\Rd) = \sigma(\Ld) = \{z\in\C,\:|z| \leq \rho(\Rd)\}.
\end{equation} 
In this special case the operators $\Rd$ and $\Ld$ are adjoint, while this property in general is not true for $\St$ and $\Vu$.\\

 The equalities (\ref{eq:1.4}), (\ref{eq:1.5}) are the analogue of (\ref{eq:1.3}) in $\lw$
however our proof is quite different from that in \cite{R1} and we use essentially Theorem 2. 
Moreover, these results agree with the spectrum of composition operator studied in \cite{R2} and the circular symmetry about 0. 
In the standard case $\omega = 1$, the spectral results (\ref{eq:1.4}), (\ref{eq:1.5}) are well known 
(see, for example Chapter V, \cite{EN}). Their proof in this  special case is based on the fact that the spectrum of 
the generator $A$ of $(\St_t)_{t\geq 0}$ is in 
$\{z\in \C, \:\Re z \leq 0\}$
 and the spectral mapping theorem for semigroups yields $\sigma(\St_t) = \{z\in \C,\:|z| \leq 1\}.$ 
Notice also that in this case we have
$$s(A)=\sup \{\Re \lambda:\: \lambda \in \sigma(A)\} = \alpha_0 = 0,$$
so the spectral bound $s(A)$ of $A$ is equal to the ground order and there is no {\it spectral gap}. 
In the general setting we deal with it is quite difficult to describe the spectrum of $A$. 
Consequently, we cannot obtain (\ref{eq:1.4}) from the spectrum of $A$ and our techniques are not based on $\sigma(A).$ 
Moreover, if for the semigroup $(\St_t)_{t\geq 0}$ on $L_{\omega}^2(\R^+)$ we can apply the spectral mapping theorem, since $\St_t$ preserves positive functions (see \cite{We1}, \cite{We2}), in general this is not true for other Hilbert spaces of functions and we could have a spectral gap $s(A) < \alpha_0$. This shows the importance of our approach which works also for more general Hilbert spaces $H$ of functions (see the conditions on $H$ listed below).
To our best knowledge it seems that Theorem 3 is the first result in the literature giving a complete characterization of 
$\sigma(\St_t)$
 and $\sigma(\Vt)$ on the spaces $L^2_{\omega}(\R^+)$. On the other hand, for the weighted two-sided shift ${\bf S}$ in $L^2_{\omega}(\R)$ a similar result has been established in \cite{V8} saying that
$$\sigma ({\bf S}) = \{ z \in \C: \: \frac{1}{\rho({\bf S}^{-1})} \leq |z| \leq \rho ({\bf S})\}.$$
 
Following the arguments in \cite{V9}, the results of this paper may be extended to a larger setup.
 Indeed, instead of $\lw$ we may consider a Hilbert space $H$ of functions on $\rr$ satisfying the following conditions:\\

(H1) $\ccp \subset H \subset \LL$, with continuous inclusions, and $\ccp$ is dense in $H$.\\

(H2) For every $t \in \R$, we have $\St_t(H)\subset H$ (resp. $\Vt(H)\subset H) $ and\\ $\sup_{t\in K}\|\St_t\|<+\infty$ (resp. $\sup_{t\in K}\|\Vt\|<+\infty$),
 for every compact set $K\subset \rr$. \\

(H3) For every $\alpha \in \R$, let $E_\alpha$ be the operator defined by 

$$E_\alpha:H\ni f\longrightarrow \Big(\R\ni x\longrightarrow f(x)e^{i\alpha x}\Big).$$

 We have $E_\alpha(H)\subset H$ and moreover,
 $\sup_{\alpha \in \R}\| E_\alpha \|<+\infty.$\\
 
 (H4) There exist $C_1 > 0$ and $a_1 \geq 0$ such that $\|R_t\|\leq C_1e^{a_1|t|},\:\forall t\in\R^+.$\\

(H5) There exist $C_2 > 0$ and $a_2 \geq 0$ such that $\|\Vt\|\leq C_2e^{a_2|t|},\:\forall t\in\R^+.$\\

Taking into account (H3), without lost of generality we may consider that in $H$ we have $\|fe^{i\alpha.}\|=\|f\|.$
For the simplicity of the exposition we deal with the case $H=\lw$ and the reader may consult \cite{V9} for the changes
necessary to cover the more general setup.\\

\section{Proof of Theorem 1}

Denote by $A$ the generator of the semi-group $(\St_t)_{t\geq 0}$.
By using the arguments based on the spectral results for semigroups (see \cite{G}, \cite{IH}) we will prove the following
\begin{lem} Let  $\lambda$ be such that $e^\lambda\in \sigma(\St)$ and $\Re \lambda=\alpha_0$. Then
there exists a sequence $(n_k)_{k \in \N}$ of integers and a sequence $(f_{m_k})_{k\in \N}$ of functions of $H$ such that  
\begin{equation}\label{eq:2.1}
\lim_{k\to\infty}\|\Bigl(e^{t A} - e^{(\lambda + 2 \pi i n_k)t}\Bigr) f_{m_k}\|=0, 
\:\forall t \in \R^+,\:\:\|f_{m_k}\|=1,\:\forall k\in \N.
\end{equation}
\end{lem}



{\bf Proof.} We have to deal with two cases: (i) $\lambda \in \sigma(A)$, (ii) $\lambda \notin \sigma (A).$
 In the case (i) $\lambda$ is in the approximative point spectrum of $A$. This follows from the fact that for any $\mu \in \C$ with
$\Re \mu > \alpha_0$ we have $\mu \notin \sigma(A)$, since $s(A) \leq \alpha_0$. 

Let $\mu_m$ be a sequence such that $\mu_m\to \lambda,\: \Re \mu_m > \alpha_0$.
 Then $\|(\mu_m I - A)^{-1} \| \geq ({\rm dist}\: (\mu_m , spec (A)))^{-1}$, hence $\|(\mu_m I- A)^{-1}\| \to  \infty.$
 Applying the uniform boundedness principle and passing to a subsequence $\mu_{m_k}$, we may find $f \in H$ such that 
$$\lim_{k \to \infty} \|(\mu_{m_k}I - A)^{-1} f\| = \infty.$$
 Introduce $f_{m_k} \in D(A)$ defined by
$$ f_{m_k} =  \frac{(\mu_{m_k}I - A)^{-1} f }{\|(\mu_{m_k}I - A)^{-1} f \|}.$$
The identity 
$$(\lambda - A) f_{m_k} = (\lambda-\mu_{m_k}) f_{m_k} + (\mu_{m_k} - A) f_{m_k}$$ 
implies that $(\lambda - A)f_{m_k} \to 0$ as $k \to \infty.$ 
Then the equality
$$(e^{tA} - e^{t\lambda})f_{m_k} = \Bigl(\int_0^t e^{\lambda(t-s)}e^{As} ds\Bigr) (A - \lambda)f_{m_k}$$
yields (\ref{eq:2.1}), where we take $n_k = 0.$\\
To deal with the case (ii), we repeat the argument in \cite{V9} and for the sake of completeness we present the details.
 We have $e^{\lambda} \in \sigma (e^A) \setminus e^{\sigma (A)}.$ Applying the results for the spectrum of a semigroup in Hilbert 
space (see \cite{G}, \cite{IH}), we conclude that there exists a sequence of integers $(n_k)$ such that $|n_k|  \to \infty$ and
  $$\|(A - (\lambda + 2 \pi i n_k)I)^{-1}\| \geq k,\: \forall k \in \N.$$
We choose a sequence $(g_{m_k}) \in H$, $\|g_{m_k} \| = 1$ so that
$$\|(A - (\lambda + 2 \pi i n_k)I)^{-1} g_{m_k}\| \geq k/2,\: \forall k \in \N$$
and define
$$f_{m_k} =  \frac{(A - (\lambda + 2 \pi i n_k)I)^{-1} g_{m_k}}{\|(A - (\lambda + 2 \pi i n_k)I)^{-1} g_{m_k}\|}.$$
Next we have
$$(e^{tA} - e^{(\lambda + 2 \pi i n_k)t}) f_{m_k} = \Bigl(\int_0^t e^{(\lambda + 2 \pi i n_k)(t-s)} e^{sA} ds\Bigr)
 (A - (2 \pi i n_k+\lambda)I)f_{m_k}$$
and we deduce (\ref{eq:2.1}).
$\Box$


\begin{lem} Let $\lambda$ be such that $e^\lambda\in \sigma(\St)$ and $\Re \lambda = \alpha_0$. Then,
there exists a sequence $(n_k)_{k \in \N}$ of integers and a sequence $(f_{m_k})_{k\in \N}$ of functions of $H$ such 
that for all $t\in \R$, 
\begin{equation}\label{eq:2.2}
\lim_{k\to \infty}\Bigl\|\Bigl(\rp S_t \ell_0- e^{(\lambda + 2 \pi i n_k)t}\Bigr) f_{m_k}\Bigr\|=0,\:\:\|f_{m_k}\|=1,\:\forall k\in \N.
\end{equation}
\end{lem}

{\bf Proof.} Clearly, for $t\geq 0$ we get (\ref{eq:2.2}) by (\ref{eq:2.1}).
Moreover, we have 
$$\|(\Vt-e^{-(\lambda + 2 \pi i n_k)t})f_{m_k}\|= \|(\Vt-e^{-(\lambda + 2 \pi i n_k)t}\Vt \St_t)f_{m_k}\|$$
$$\leq \|\Vt\||e^{-(\lambda + 2 \pi i n_k)t}|\Big\|\Bigl(e^{(\lambda + 2 \pi i n_k)t}-\St _t \Bigr)f_{m_k}\Big\|,
\:\forall t\in \rr.$$
Thus
$$\lim_{k\to \infty} \|(\Vt-e^{-(\lambda + 2 \pi i n_k)t})f_{m_k}\|=0$$
and this completes the proof of (\ref{eq:2.2}).
$\Box$\\
Recall that for $\phi\in \Cc$, $T_\phi$ is the operator on $\lw$ given by 
$$T_\phi(f)=\rp (\phi*\ell_0 f), \:\forall f\in \lw.$$

\begin{lem}
For all $\phi\in \Cc$ and $\lambda$ such that $e^\lambda \in \sigma (\St )$ with $\Re \lambda = \alpha_0$ we have 

\begin{equation}{\label{eq:eqi2}}
|\hat{\phi}(i\lambda +a)|\leq \|T_\phi\|,\:\forall a\in \R.
\end{equation}

\end{lem}

{\bf Proof.} Let $\lambda\in \C$ be such that $e^\lambda\in \sigma (\St )$ with  $\Re \lambda = \alpha_0$ and let $(f_{m_k})_{k \in \N}$ 
be the sequence satisfying (\ref{eq:2.2}). Fix $\phi \in \Cc$ and consider

$$|\hat{\phi}(i\lambda + a)| =  |\int_{\R} \langle \phi(t) e^{(\lambda -i a)t}  f_{m_k}, f_{m_k} \rangle dt|$$
$$\leq \Big |\int_{\R} \big \langle \phi(t)\Bigl( e^{(\lambda + i 2 \pi n_k  )t} - \rp S_t \ell_0\Bigr)
 e^{-i(a + 2 \pi n_k)t} f_{m_k}, f_{m_k} \big \rangle dt\Big| $$
$$+ \Big | \int_{\R} \langle \phi(t)\rp  S_t \ell_0  e^{-i(a + 2 \pi n_k)t} f_{m_k}, f_{m_k} \rangle dt\Big |.$$

The first term on the right side of the last inequality goes to 0 as $k \to \infty$ since by Lemma 1, for every fixed $t$ we have
$$\lim_{k\to +\infty}\big \Vert e^{-i(a + 2 \pi n_k)t}\Bigr(e^{(\lambda + i 2 \pi n_k )t} - \rp S_t \ell_0 \Bigr)f_{m_k} \big \Vert=0.$$
On the other hand, 
$$I_{k} =\Big |\int_\R <\phi(t) \rp S_t \ell _0  e^{-i(a + 2 \pi n_k)t} f_{m_k}, f_{m_k}> dt \Big|$$
$$ = \Big|\langle \rp \Bigl[ \int_\R \phi(t) e^{-i(a + 2 \pi n_k)t} \ell_0 f_{m_k}(. - t)dt\Bigr], f_{m_k}(.) \rangle\Big|$$
$$= \Big|\langle \rp  \int_\R  \phi(. - y) e^{i( a +2 \pi n_k) y}\ell_0 f_{m_k}(y) dy, e^{ i( a + 2 \pi n_k) .}f_{m_k}(.)\rangle\Big|$$
$$ = \Big|\langle \Bigl(T_{\phi} (e^{i(a + 2 \pi n_k).} f_{m_k} )\Bigr), e^{ i(a + 2 \pi n_k).} f_{m_k}(.)\rangle\Big|$$
and $I_{k} \leq \|T_{\phi}\|.$ Consequently, we deduce that
$|\hat{\phi}(i\lambda + a)| \leq \|T_{\phi}\|. \:\Box$

Notice that the property (\ref{eq:eqi2}) implies that 
$$|\hat{\phi}(\lambda)|\leq \|T_\phi\|,\:\forall \lambda\in \C,\:{\rm provided} \:\Im \lambda=\alpha_0.$$

\begin{lem}
Let $\phi\in \Cc$ and let $\lambda$ be such that $e^{-\bar{\lambda}} \in \sigma ((\Vu)^*)$ with $ \
Re\: \lambda = -\alpha_1$. Then  we have
\begin{equation}\label{eq:eqia2}
|\hat{\phi}(i\lambda +a)|\leq \|(T_\phi)\|,\:\forall a\in \R.
\end{equation}

\end{lem}

{\bf Proof.} Consider the semigroup $(\Vt)^*_{t\geq 0}$ and  let $B$ be its generator. 
We identify $H$ and its dual space $H'$. So the semigroup $(\Vt)^*$, $t\geq 0$ is acting on $H$.
Let $\lambda\in \C$ 
be such that $e^{-\bar{\lambda}}\in \sigma ((\Vu)^*)$ and $|e^{-\bar{\lambda}}|=\rho(\Vu) = \rho ((\Vu)^*) = e^{\alpha_1}.$
Then, by the same argument as in Lemma 1, we prove that
 there exists a sequence $(n_k)_{k \in \N}$ of integers and a sequence $(f_{m_k})_{k\in \N}$ of functions of $H$
 such that for all $t\in \rr$,
$$\lim_{k\to\infty}\|(e^{tB}-e^{(-\bar{\lambda} + i 2 \pi  n_k)t})f_{m_k}\|=0$$
and $\|f_{m_k}\|=1$.
Thus we deduce
$$\lim_{ k\to +\infty}\|(\Vt)^* f_{m_k}-e^{-(\bar{\lambda}-i 2\pi n_k)t} f_{m_k}\|=0,\: t \geq 0.$$
Since for $ t \geq 0$ we have $\Vt \St _t =I$, we get $(\St _t  )^*(\Vt)^*=I$.
Then, for $t\geq 0$ we get
$$\|(\St _t )^* f_{m_k}-e^{(\bar{\lambda}-i 2\pi n_k)t} f_{m_k}\|$$
$$=\|(\St _t )^* f_{m_k}-e^{(\bar{\lambda}-i2 \pi n_k)t}(\St _t )^*(\Vt)^* f_{m_k}\|$$
$$\leq\|(\St _t )^*\||e^{(\bar{\lambda}-i 2\pi n_k)t}|\|(e^{-(\bar{\lambda}-i 2\pi n_k)t}f_{m_k}-(\Vt)^*f_{m_k})\|.$$
This implies that 
\begin{equation} \label{eq:2.5}
\lim_{k\to+\infty}\|((\rp  S_t \ell_0)^* -e^{(\bar{\lambda}-i 2\pi n_k)t}) f_{m_k}\|=0,\:\forall t\in \R.
\end{equation}
 We write
$$\hat{\phi}(i\lambda + a) = \int_{\R} < \phi(t) e^{- i (a +2 \pi  n_k) t} f_{m_k}, 
e^{\bar{\lambda} t - 2 i\pi  n_k t} f_{m_k} >dt$$
$$ = \int_{\R} < \phi(t) e^{-i( a + 2 \pi  n_k) t} f_{m_k}, \Bigr(e^{(\bar{\lambda}  - i 2 \pi n_k) t} - (\rp  S_t \ell_0)^* \Bigr)f_{m_k} > dt$$
$$+ \int_{\R} < \phi(t)e^{ -i(a + 2 \pi  n_k) t}(\rp  S_t \ell_0)f_{m_k}, f_{m_k} > dt = J_{k}' + I_{k}'.$$
 From (\ref{eq:2.5}) we deduce that $J_{k}' \to 0$ as $k \to \infty.$ 
For $I_{k}'$ we apply the same argument as in the proof of Lemma 3 and we get
$|\hat{\phi}(i\lambda)| \leq \|T_{\phi}\|.\:\:\Box$

\vspace{0.2cm}

\begin{lem}
For every function $\phi \in \Cc$ and for $z \in U=\{z\in \C, \:\Im z\in [-\alpha_1, \alpha_0]\}$ we have
$$|\hat{\phi}(z)|\leq\|T_\phi\|.$$ 
\end{lem}
{\bf Proof.} We will use the Phragm\'en-Lindel\"of theorem and we start by proving the estimations on the bounding lines. 
There exists $\alpha = e^{-i z}\in \sigma (\St )$ such that $|\alpha|= e^{\Im z} = e^{\alpha_0}$. Following (\ref{eq:eqi2}), we obtain 
$$|\hat{\phi}(z)|\leq \|T_\phi\|,$$
for every $z$ such that $\Im z = \alpha_0$. Next notice that $\rho(\Vu) = \rho \Bigl((\Vu)^*\Bigr).$ So there exists $\beta = e^{-i \bar{z}} = e^{-\overline{(-iz)}} \in \sigma ((\Vu)^*)$ such that $|\beta|= e^{\alpha_1}$ and
$$- \Im z = \ln |\beta| = \alpha_1.$$
Then taking into account (\ref{eq:eqia2}), we get
$$|\hat{\phi}(z)|\leq \|T_\phi\|,$$
for every $z$ such that $\Im z= - \alpha_1$. In the case $\alpha_1 + \alpha_0 = 0$ the result is obvious. 
So assume that $\alpha_0 + \alpha_1 > 0.$
 Since $\phi\in \Cc$ we have 
$$|\hat{\phi}(z)|\leq C\|\phi\|_{\infty}e^{k |\Im z|}\leq K \|\phi\|_{\infty},\:\:\forall z \in U,$$
where $C > 0$, $k >0$ and $K >0$ are constants.
An application of the Phragm\'en-Lindel\"of theorem for the holomorphic function $\widehat{\phi}$, yields
$$|\widehat{\phi}(\alpha)|\leq \|T_\phi\|$$
for $\alpha\in \{ z\in \C: \Im z \in \I\}.$
$\Box$\\

Notice that for multipliers in $L^p(\R)$ studied in \cite{H} some  relations concerning the norm $\|T_{\phi}\|$ and $\phi$ hold.
Combining the results in Lemma 3-5, we get
\begin{lem}
For every $\phi \in C_c^\infty(\R)$ and for every $a\in \I$  we have
 $$|\widehat{(\phi)_a}(x)|\leq \|T_\phi\|,\:\forall x\in \R.$$
 
\end{lem}





{\bf Proof of Theorem 1.} 
   Let $T$ be a Wiener-Hopf operator. Let $a \in \I$. The proof follows the approach in \cite{V2}. Notice that in \cite{V2}, we establish (\ref{eq:r})
 for $f\in C_c^\infty(\rr)$ and which is new here is that we prove (\ref{eq:r}) for $f\in \lw$ such that $(f)_a\in L^2(\rr)$. 
Following \cite{V2}, there exists a sequence $(\phi_n)_{n \in \N}\subset \Cc$ such
 that $T$ is the limit of $(T_{\phi_n})_{n \in \N}$ with respect to the strong operator topology
and we have  $\|T_{\phi_n}\|\leq C \|T\|,$ where $C$ is a constant independent of $n$. 
 According to Lemma 6, we have
\begin{equation}\label{eq:comp}
|\widehat{(\phi_n)_a}(x)|\:\leq \| T_{\phi_n}\| \leq C\|T\|,\:\forall x \in \R,\:\forall n \in \N
\end{equation}
 and we replace $(\widehat{(\phi_n)_a})_{n \in \N}$ by a suitable subsequence, also denoted by $(\widehat{(\phi_n)_a})_{n \in \N}$, 
converging with respect to the weak topology $\sigma (L^{\infty}(\R),L^1(\R))$ to a function $h_a\in L^{\infty}(\R)$ such that
$\|h_a\|_\infty\leq \:C\:\|T\|.$
We have
$$\lim_{n \to +\infty} \int_\R{\Bigl(\widehat{(\phi_n)_a}(x)-h_{a}(x)\Bigr)\:g(x)\:dx}=0, \:\:\forall g \in L^1(\R).$$
Fix $f\in \lw$ so that $(f)_a\in L^2(\R^+).$
Then we get
$$\lim_{n \to +\infty} \int_{{\R}}{\Bigl(\widehat{(\phi_n)_a}(x)\widetilde{(f)_a}(x)-h_{a}(x) \widetilde{(f)_a}(x)\Bigr)\:g(x)\:dx}=0,$$
 for all  $g \in L^2({\R})$. We conclude that $\Bigl(\widehat{(\phi_n)_a}\widetilde{(f)_a}\Bigl)_{n \in\N}$ converges weakly in ${L^2({\R})}$ to $h_{a}\widetilde{(f)_a}$.\\
 On the other hand, we have
 $$(T_{\phi_n}f)_a=\rpb  ((\phi_n)_a*\ell_0 (f)_a)=\rpb  \F^{-1}(\widehat{(\phi_n)_a}\widetilde{(f)_a})$$
and thus $(T_{\phi_n}f)_a$ converges weakly in $L^2(\rr)$ to $\rpb  \F^{-1}(h_a \widetilde{(f)_a})$. 
For $g\in \Cc$, we obtain 
$$\int_{\R^+}\Bigl|(T_{\phi_n}f)_a(x)-(Tf)_a(x)\Bigr|\:|g(x)|\:dx $$
$$\leq C_{a,g} \| T_{\phi_n}f-Tf \|,\:\forall n \in \N,$$
where $C_{a,g}$ is a constant depending only of $g$ and $a$.
Since $(T_{\phi_n}f)_{n \in \N}$ converges to $Tf$ in ${L_{\omega}^2({\R^+})}$, we get
$$\lim_{n \to +\infty}\int_{\R^+}{(T_{\phi_n}\:f)_a(x)g(x)\:dx}= \int_{\R^+}{(Tf)_a(x)g(x)\:dx},\:\:\forall g \in \Cc.$$
Thus we deduce that $(Tf)_a=\rpb  \F^{-1}(h_a \widetilde{(f)_a})$. 
The symbol $h$ is holomorphic on $\overset{\circ}{U}$ following the same arguments as in \cite{V2}. 
$\Box$

\section{Preliminary spectral result}

As a first step to our spectral analysis in this section we prove the following

\begin{prop}
Let $T\in \M$ and suppose that the symbol $\mu$ of $T$ is continuous on $U$. 
Then ${\mu(U)} \subset \sigma(T).$
\end{prop}

{\bf Proof of Proposition 1.} 
Let $T$ be a bounded operator on $H$ commuting with $\St _t,\: t \geq 0$ or $\Vt, \: t \geq 0$. For $a\in \I$, we have
$$(Tf)_a=\rpb  \F^{-1}(\mu_a \widetilde{(f)_a}),\:\forall f\in \lw,$$
where $\mu_a \in L^{\infty}(\R)$,
provided $(f)_a\in L^2(\rr)$.
Suppose that  $\lambda \notin \sigma (T)$. 
Then, it follows easily that the resolvent $(T - \lambda I)^{-1}$ also commutes with $(\St _t)_{t\in\rr}$ or $(\Vt)_{t\in \rr}$. 
Consequently,
$(T-\lambda I)^{-1}$
 is a Wiener-Hopf operator and  for $a\in \I$ there exists a function $h_a \in L^{\infty}(\R)$ such that
$$((T-\lambda I)^{-1}g)_a=\rpb  \F^{-1}(h_a\widetilde{(g)_a}),$$
for $g\in \lw$ such that $(g)_a\in L^2(\rr).$ 
If $f$ is such that $(f)_a\in L^2(\R^+)$, set $g=(T-\lambda I)f$. Then following Theorem 1, we deduce that $(Tf)_a \in L^2(\R^+)$ and $(g)_a=((T-\lambda I)f)_a \in L^2(\rr)$. 
Thus applying once more Theorem 1, we get
$$((T-\lambda I)^{-1}(T-\lambda I)f)_a=\rpb  \F^{-1}(h_a\widetilde{(T-\lambda I)f)_a}$$
$$=\rpb  \F^{-1}\Bigl(h_a\F \rpb  [\F^{-1}((\mu_a-\lambda) \widetilde{(f)_a})]\Bigr).$$

We have
$$\|(f)_a\|_{L^2}\leq \|h_a\F \rpb  \F^{-1}((\mu_a-\lambda) \widetilde{(f)_a})\|_{L^2}
\leq \|h_a\|_\infty \|\F \rpb  \F^{-1}((\mu_a-\lambda) \widetilde{(f)_a})\|_{L^2}$$
and we deduce 
\begin{equation}\label{eq:contra}\|\widetilde{(f)_a}\|_{L^2}\leq C \|(\mu_a-\lambda) \widetilde{(f)_a}\|_{L^2},
\end{equation}
for all $f\in \lw$ such that $(f)_a\in L^2(\rr)$. 
 Let $\lambda=\mu_a(\eta_0) = \mu(\eta_0 + i a) \in \mu (U)$ for $a \in [- \ln \rho(\Vu), \: \ln \rho(\St )]$ and some $\eta_0\in \R$. Since
 the symbol $\mu$ of $T$ is continuous, the function $\mu_a(\eta) = \mu(\eta + i a)$ is continuous on $\R.$ We will construct
a function $f(x) = F(x) e^{-ax}$ with ${\rm supp} (F)\: \subset \R^+$ for which (\ref{eq:contra}) is not fulfilled.
Consider 
$$g(t) = e^{-\frac{b^2(t - t_0)^2}{2}} e^{i(t-t_0)\eta_0},\: b > 0, t_0 > 1$$
with Fourier transform
$$\hat{g}(\xi) = \frac{1}{b} e^{-\frac{(\xi - \eta_0)^2}{2b^2}} e^{-it_0 \xi}.$$

Fix a small $0 < \epsilon < \frac{1}{2}C^{-2}$, where $C$ is the constant in (\ref{eq:contra}) and let 
$\delta > 0$ be fixed so that $|\mu_a(\xi) - \lambda| \leq \sqrt{\epsilon}$ for $\xi \in V = \{\xi \in \R:\: |\xi - \eta_0| \leq \delta\}.$ 
Moreover, assume that
$$|\mu_a(\xi) - \lambda|^2 \leq C_1,\: a.e.\: \xi \in \R.$$
We have for $0 < b \leq 1$ small enough
$$\int_{\R \setminus V} |\hat{g}(\xi)|^2 d\xi \leq \frac{1}{b^2} \int_{|\xi - \eta_0| \geq \delta} e^{-\frac{(\xi - \eta_0)^2}{2b^2}} d\xi$$
$$\leq e^{-\frac{\delta^2}{4b^2}} \frac{1}{b^2} \int_{|\xi - \eta_0| \geq \delta} e^{-\frac{(\xi - \eta_0)^2}{4b^2}} d\xi \leq C_0 b^{-1} 
e^{-\frac{\delta^2}{4b^2}} \leq \epsilon$$
with $C_0 > 0$ independent of $b > 0$.
 We fix $b > 0$ with the above property and we choose a function $\varphi \in C_c^{\infty}(\R^+)$ such that
$0 \leq \varphi \leq 1, \: \varphi(t) = 1$ for $1 \leq t \leq 2t_0 - 1,\:\varphi(t) = 0$ for $t \leq 1/2$ and for
$t \geq 2t_0 - 1/2.$ We suppose that $|\varphi^{(k)}(t)|\leq c_1,\: k = 1,2, \forall t \in \R$.
Set $G(t) = (\varphi(t) - 1) g(t).$ We will show that
\begin{equation} \label{eq:outi}
|(1 + \xi^2) \hat{G}(\xi)| \leq\sqrt{\frac{C_2}{4\pi} \epsilon}
\end{equation}
for $t_0$ large enough with $C_2 > 0$ independent of $t_0.$  On the support of $(\varphi - 1)$ we have $|t - t_0| > t_0 -1$ and integrating by parts in $\int_{\R} (1 + \xi^2) G(t) e^{-it\xi} dt$ we must estimate the integral
$$\int_{|t - t_0| \geq t_0 -1} e^{-\frac{b^2(t-t_0)^2}{2}} (1 + |t - t_0| + (t - t_0)^2)dt $$
$$\leq \Bigr(\int_{-\infty}^{1 - t_0}(1+ |y| + y^2)e^{-b^2 y^2/2} dy + \int_{t_0 -1}^{\infty} (1 + y + y^2) e^{-b^2 y^2/2} dy\Bigr).$$
Choosing $t_0$ large enough we arrange (\ref{eq:outi}).\\

We set $F = \varphi g \in C_c^{\infty}(\R^+)$ and we obtain
$$\int_{\R \setminus V} |\tilde{F}(\xi)|^2 d\xi \leq 2\int_{\R \setminus V} |\hat{g}(\xi)|^2 d\xi + 2 \int_{\R \setminus V} |\hat{G}(\xi)|^2 d\xi$$
$$ \leq 2\epsilon + \frac{C_2 \epsilon}{2\pi} \int_{\R} (1 + \xi^2)^{-2} d\xi \leq (2
 + C_2)\epsilon.$$ 

 Then
$$\int_{\R} |(\mu_a(\xi) - \lambda)\tilde{F}(\xi)|^2 d\xi \leq \int_{\R \setminus V} |(\mu_a(\xi) - \lambda)\tilde{F}(\xi)|^2 d\xi+ \int_V |(\mu_a(\xi) - \lambda)\tilde{F}(\xi)|^2 d\xi $$
$$\leq C_1 (2 + C_2)\epsilon + (2\pi)^2\|F\|^2_{L^2}\epsilon.$$
Now assume (\ref{eq:contra}) fulfilled. Therefore 
$$ (2 \pi)^2\|F\|_{L^2}^2 \leq C^2\|(\mu_a(\xi) - \lambda)\hat{F}(\xi)\|_{L^2}^2 \leq C^2C_1(2 + C_2) \epsilon + (2 \pi C)^2 \|F\|^2_{L^2} \epsilon,$$
and since $C^2\epsilon < \frac{1}{2}$, we conclude that
$$\|F\|_{L^2}^2 \leq \frac{C^2C_1}{2 \pi^2}(2 + C_2) \epsilon .$$

On the other hand,
$$\|F\|_{L^2}^2 \geq \frac{1}{2}\|g\|_{L^2}^2 - \|(\varphi -1) g\|_{L^2}^2 \geq \frac{1}{2} \|g\|_{L^2}^2 - (2\pi)^{-2} \frac{C_2}{2} \epsilon$$
and
$$\int_{\R} |g(t)|^2 dt \geq \int_{|t - t_0| \leq \frac{1}{b}} e^{-b^2(t- t_0)^2}dt \geq  \frac{2e^{-1}}{b} \geq 2e^{-1}.$$
 For small $\epsilon$ we obtain a contradiction, since $C_2$ is independent of $\epsilon$. This completes the proof. $\Box$



\section{Spectra of $(\St _t)_{t\in \rr}$, $(\Vt)_{t\in \rr}$ and bounded operators commuting with at least one of these semigroups}

Observing that the symbol of $\St _t$ is 
$z \longrightarrow e^{-i tz}$, an application of Proposition 1 to the operator $\St _t$ yields
\begin{equation} \label{eq:4.1}
\{z\in \C,\:e^{-\alpha_1 t}\leq |z|\leq e^{\alpha_0t}\}\subset \sigma (\St _t).
\end{equation}
This inclusion describes only a part of the spectrum of $\St _t$. We will show that  in our general setting we have (\ref{eq:1.4}). 
To prove this, for $t>0$  assume that $z\in \C$ is such that
$ 0< |z| < e^{-\alpha_1 t}.$
Let $g\in H$ be a function such that $g(x) = 0$ for $x \geq t$ and  $g \not= 0$.
 If the operator $(zI - \St _t)$ is surjective on $H$, then there exists
$f \not= 0$ such that $(z- \St _t)f = g$. This implies
$\Vt g = 0$ and hence
$$\Bigl(\Vt - \frac{1}{z} I\Bigr) f = 0$$
which is a contradiction. So every such $z$ is in the spectrum of $\St _t$ and we obtain (\ref{eq:1.4}).

Next, it is easy to see that in our setup for the approximative point spectrum $\Pi(\St _t)$ of $\St _t$ 
we have the inclusion
\begin{equation} \label {eq:4.2}
\Pi(\St _t) \subset \{z\in \C:\: e^{-\alpha_1 t} \leq |z|\leq e^{\alpha_0 t}\}.
\end{equation}
Indeed, for  $z \neq 0,$ we have the equality
$$\St _{-t} - \frac{1}{z}I = \frac{1}{z} \St _{-t} (z I - \St _t).$$

If for $z \in \C$ with $0 < |z| < e^{-\alpha_1 t},$ there exists a sequence $(f_n)$ such that $\|f_n\| =1 $ 
and $\|(z I - \St_t)f_n\| \to 0$ as $n \to \infty$, then
$$\Bigl(\Vt - \frac{1}{z} I\Bigr) f_n \to 0,\: n \to \infty$$
and this leads to $\frac{1}{z} \in \sigma (\Vt)$ which is a contradiction. Next, if $0 \in \Pi(\St _t)$, 
there exists a sequence $g_n \in H$ such that $\St_t g_n \to 0,\: \|g_n\| = 1.$ 
Then $g_n = \Vt \St _t g_n$ and we obtain a contradiction. Since the symbol of $\Vt$ is $z \longrightarrow e^{i t z}$, applying Proposition 1, we obtain
$$\{z\in \C:\: e^{-\alpha_0 t} \leq |z|\leq e^{\alpha_1 t} \}\subset \sigma (\Vt).$$
Passing to the proof of (\ref{eq:1.5}), notice that $\St _t^*(\Vt)^* = I.$ Then for $0 < |z| < e^{-\alpha_0 t}$ we have
\begin{equation} \label{eq:4.3}
z \Bigl(\frac{1}{z}I - \St _t^*\Bigr) = \St _t^*\Bigl( (\Vt)^* - z I \Bigr).
\end{equation}
It is clear that $0 \in \sigma_r(\St _t)$, where $\sigma_r(\St _t)$ 
denotes the residual spectrum of $\St _t$. In fact, if $0 \notin \sigma_r(\St_t)$, 
then 0 is in the approximative point spectrum of $\St_t$ and this contradicts (\ref{eq:4.2}). 
Since $0 \in \sigma_r(\St_t)$, we deduce that 0 is an eigenvalue of $\St_t^*$. 
Let $\St _t^* g = 0, \: g \not= 0.$ Assume that $(\Vt)^* - z I$ is surjective.
 Therefore, there exists $f \not= 0$ so that $((\Vt)^* - z) f = g$ 
and (\ref{eq:4.3}) yields $(\frac{1}{z}- \St _t^*) f = 0.$
Consequently, $\frac{1}{|z|} \leq \rho(\St_t^*) = \rho(\St _t) = e^{\alpha_0 t}$ 
and we obtain a contradiction. Thus we conclude that $z \in \sigma ((\Vt)^*)$, hence $\bar{z} \in \sigma (\Vt)$ 
and the proof of (\ref{eq:1.5}) is complete.\\

To study the operators commuting with $(\St _t)_{t\in \rr}$, we need the following


\begin{lem}
 Let $\phi\in \CC$. The operator $T_\phi$ commutes with $\St_t, \: \forall t > 0,$ if and only if the support of $\phi$ 
is in $\overline{\R^+}$. 
\end{lem}
{\bf Proof.} First if $\psi \in L^2_{\omega}(\R^+)$ has compact support in $\overline{\R^+}$, it is easy to see that $T_{\psi}$ commutes with $\St _t, \: t \geq 0$. Now consider $\phi\in \CC$ and suppose that $T_{\phi}$ commutes with $\St _t,\: t \geq 0.$ We write 
$\phi = \phi \chi_{\R^{-}} + \phi \chi_{\R^{+}}.$ If $T_{\phi}$ commutes with $\St _t,\: t \geq 0$, then the operator $T_{\phi \chi_{\R^{-}}}$ commutes too.
 Let the function $\psi = \phi \chi_{\R^{-}}$ have support  in $[-a, 0]$ with $a > 0$.  
Setting $f=\chi_{[0, a]}$, we get $S_af=\chi_{[a, 2a]}$. 
For $x\geq 0$ we have 
$$T_\psi R_a r_+f(x)=\rpb  (\psi*S_a f) (x)= \int_{-a}^0 \psi(t) \chi_{\{a\leq x-t\leq 2a\}}dt
=\int_{max (-a, -2a+x)}^{min(x-a, 0)}\psi(t) dt.$$
Since $T_\psi \St_a=\St_a T_\psi$, we deduce
$\int_{-a}^{x-a} \psi(t) dt=0,\:\forall x\in [0, a].$
This implies that 
$\psi(t)=0,$ for $t\in [-a, 0]$ and $supp(\phi)\subset \overline{\R^+}.$ 
$\Box$

\begin{lem} Let  $\lambda$ be such that $e^\lambda\in \sigma (\St )$. Then
there exists a sequence $(n_k)_{k \in \N}$ of integers and a sequence $(f_{m_k})_{k\in \N}$ of functions of $H$ such that  
\begin{equation}\label{eq:2.1bis}
\lim_{k\to\infty}<\Bigl(\St _t - e^{(\lambda + 2 \pi i n_k)t}\Bigr) f_{m_k}, f_{m_k}>=0, 
\:\forall t \in \R^+,\:\:\|f_{m_k}\|=1,\:\forall k\in \N.
\end{equation}
\end{lem}
{\bf Proof.} If $\lambda \notin \sigma(A)$, we repeat the argument of the proof of Lemma 1. 
Denote by $\sigma_r(A)$ the residual spectrum  of $A$.
If $\lambda \in \sigma(A) \setminus  \sigma_r(A)$, we deduce that $\lambda$ is in the approximative point spectrum of $A$ and we can apply the argument of Lemma 1. Finally, if $\lambda \in \sigma_r(A)$, then there exists $f\in H$ such that $A^*f=\overline{\lambda} f$ and $\|f\|=1$. 
We set $f_{m_k}=f,\: n_k=0$, for $k\in \N$ and we use the fact that $(R_t)^*f = e^{\bar{\lambda} t} f.$
$\Box$


\begin{lem}
For all $\phi\in C_c^{\infty}(\overline{\R^+})$ and $\lambda$ such that $e^\lambda \in \sigma(\St )$ we have 
\begin{equation}{\label{eq:eqi1bis}}
|\hat{\phi}(i\lambda)|\leq \|T_\phi\|.
\end{equation}
\end{lem}

The proof is based on the equality
$$ \hat{\phi}(i \lambda) = \int_{\R^+} \phi(t) e^{\lambda t} dt =  \int_{\R^+} \langle \phi(t) e^{(\lambda  +  2 \pi i n_k)t}f_{m_k}, e^{2 \pi i n_k t} f_{m_k} \rangle dt$$
$$= \int_{\R^+} \langle \phi(t) \Bigl( e^{(\lambda + 2\pi i n_k )t}I - \St _t\Bigr) f_{m_k},  e^{2 \pi i n_k t} f_{m_k}\rangle dt +
\int_{\R^+}\langle  \phi(t) \St_t f_{m_k}, e^{2 \pi i n_k t} f_{m_k}\rangle dt.$$
We apply Lemma 8 and we repeat the argument of the proof of Lemma 3. Notice that here the integration is over $\R^+$ and we do not need to examine the integral for $t < 0.$\\
Following \cite{V2}, the operator $T$  is a limit of a sequence of operators $T_{\phi_n}$,
 where $\phi_n\in \CC$ and $\|T_{\phi_n}\|\leq C\|T\|.$
The sequence $(T_{\phi_n})_{n\geq 0}$ has been constructed in \cite{V2} and it follows from its construction that if 
$T$ commutes with $\St_t,\: t > 0$, then  $T_{\phi_n}$ has the same property for all $n\in \N$. 
Therefore, Lemma 7 implies that $\phi_n \in C_c^{\infty}(\overline{\R^+})$ and to obtain Theorem 2 for bounded operators commuting with 
$(\St _t)_{t>0}$, we apply Lemma 9
 and the same arguments as in the proof of Theorem 1. Finally, applying Theorem 2 and the arguments of the proof of Proposition 1,
 we establish (\ref{eq:1.6}) and this  completes the proof of iii) in Theorem 3.
Next we prove the following
\begin{lem}
Let $\phi\in \Cc$. Then $T_\phi$ commutes with $\Vt$, $\forall t>0$ if and only if 
$supp(\phi)\subset \overline{\R^-}$. 
\end{lem}
The proof of Lemma 10 is essentially the same as that of Lemma 7. By using Lemma 10, we obtain an analogue of Lemma 9 and Theorem 2 for bounded operators commuting with $(\Vt)_{t>0}$  and applying these results we establish iv) in Theorem 3.\\

{\bf Acknowledgments.} The author would like to thank the referee for helpful comments and suggestions, making the paper more understandable.



\end{document}